\documentclass[a4paper,11pt,reqno]{amsart}

\usepackage{amsmath,amssymb}
\usepackage{enumerate}
\usepackage{ifthen}
\usepackage{calc}
\nonstopmode\numberwithin{equation}{section}
\usepackage{tikz}
\theoremstyle{plain}

\theoremstyle{definition}
\newtheorem{defn}{Definition}

 \usepackage[font={footnotesize},labelfont={bf}]{caption}

\begin{document}
\title{Analogy between Learning With Error Problem and Ill-Posed Inverse Problems}
\maketitle
 \begin{center}
\center{  $\text{Gaurav Mittal}$}
\medskip
{\footnotesize

 \center{Defence Research and Development Organization, Near Metcalfe House, New Delhi, 110054, India  \\ email:  gaurav.mittaltwins@gmail.com, gmittal93.hqr@gov.in }
 }

\bigskip
\begin{abstract} In this work, we unveil an analogy between well-known lattice based learning with error problem and ill-posed inverse problems. We show that  LWE problem is a structured inverse problem. Further, we propose a symmetric encryption scheme based on ill-posed problems and thoroughly discuss its security. Finally, we propose a public key encryption scheme based on  our symmetric encryption scheme and CRYSTALS-Kyber KEM (key encapsulation mechanism)    and discuss its security.
\end{abstract}
\hspace{-63mm}\textbf{2010 MSC:} Primary: 94A60, 65J22\\
\hspace{-08mm}{\textbf{Keywords:} Learning With Error Problem, Cryptography, Inverse Problems }
\end{center}
\section{Introduction and Analogy}
\subsection{Learning with Error Problem}\ \\
\noindent The LWE (Learning With Error) problem is a well-known hard problem in cryptography \cite{Bern}. Having been the subject of intense scrutiny for nearly a decade, the crypto-primitives CRYSTALS Kyber \cite{Kyber} and CRYSTALS Dilithium \cite{Dili}, whose security inherently rely on LWE problem, have been standardized by NIST. The LWE problem was put forward by Regev in his seminal work \cite{Regev}. Formally, one may define LWE problem as follows:
\begin{defn}
For a prime $q$, let $\mathcal{D}$ be an error distribution over the modular ring $\mathbb{Z}_q$. For a given $m$ (also known as dimension parameter), let $s\leftarrow \mathbb{Z}_q^m$ be chosen uniformly at random and it is kept as a secret. For $1\leq j\leq n$, where $n$ is some polynomial in $m$, consider the samples$$
(a_j, b_j = \langle a_j, s \rangle + e_j) \in \mathbb{Z}_q^m \times \mathbb{Z}_q,
$$
where  $a_j \leftarrow \mathbb{Z}_q^n \) is picked uniformly at random. Further, $e_j\leftarrow \mathbb{Z}_q$ is picked following distribution $\mathcal{D}$ and typically, this is a short element (in terms of norm). The LWE problem is to derive $s$ from the knowledge of $(a_j, b_j)_{1\leq j\leq n}$.\\
The LWE problem can also be formulated in terms of matrix form as follows: Given
\[
\mathbf{A}= \begin{bmatrix}
(a_1)\\ (a_2)\\ \vdots\\ (a_n)
\end{bmatrix}_{m \times n}, \quad \mathbf{e}=\begin{bmatrix}
e_1\\e_2\\ \vdots\\ e_n
\end{bmatrix} \in \mathbb{Z}_q^n, \quad \mathbf{b} = \mathbf{A}\mathbf{s} + \mathbf{e}= \begin{bmatrix}
(a_1)\\ (a_2)\\ \vdots\\ (a_n)
\end{bmatrix}\times \mathbf{s} + \begin{bmatrix}
e_1\\e_2\\ \vdots\\ e_n
\end{bmatrix} \in \mathbb{Z}_q^n,
\]
we need to find $\mathbf{s}\in \mathbb{Z}_q^m$ from the knowledge of $\mathbf{A}, \mathbf{b}$.
\end{defn}\noindent 
It is well-known that LWE problem is very hard and solving it, on an average, is equivalent to solve certain lattice problems in the worst-case  using quantum computers \cite{Regev}. This makes LWE a very strong candidate for post-quantum cryptography and the worst-case advantage (discussed above) is not associated with any other post-quantum candidates \cite{Bern}. Further, to explore the roots of computational hardness of LWE problem, we switch into the subject of inverse problems. We refer to \cite{Engl} for more preliminary details on this subject. 
\subsection{Well-Posed and Ill-Posed Problems}\ \\
\noindent We consider the operator equations of the form 
\begin{equation}
\label{1}
 T:D(T)\subset U\to V \ \ \ \text{defined as}\ \ \ \ T(u)=v,\end{equation}
where $U$ and $V$ denote Banach spaces and $D(T)$ denotes the domain of $T$. Given $u$, deducing $v=T(u)$ is the direct problem. The corresponding inverse problem is deducing $u$ from the value of $v$. As per Hadamard criteria, a problem of the form (\ref{1}) is well-posed if the following three conditions are met. \begin{enumerate}[(i)]
\item For a given $v\in V$, (\ref{1}) has a solution.
\item The solution is unique.
\item Continuous dependency must be there, i.e., if $v'$ is given in place of $v$ and it is close to $v$, then the corresponding solutions $u'$ and $u$ should be near to each other (in terms of norm).
\end{enumerate}
Further, if any one of the above three conditions is not satisfied, then the  problem   (\ref{1}) is ill-posed. It is worth to note that the conditions (i)-(ii) are straight-forward to understand. The condition (iii) is less clear. Before going further, we discuss an example of an ill-posed problem, which is ill-posed since it does not fulfill condition (iii).
 \\
 \textbf{Example of an ill-posed Problem}.

\noindent Let $L^2[0,1]$ be the space of real-valued Lebesgue integrable functions $g$ defined on $[0,1]$ such that $\int_0^1 g^2(y)\, dy<\infty$. We consider an operator $$\mathcal{S} :L^2[0,1]\to L^2[0,1]$$ defined as 
 \begin{equation}
 \label{1.2}
(\mathcal{S} (\Psi))(y) = \int_0^1 e^{-|y - s|} \Psi(s) \, ds.
 \end{equation} This operator is also known as Hilbert-Schmidt operator (HSO). For the integral equation (\ref{1.2}), the kernel  $\mathcal{K}(y,s)= e^{-|y-s|}$ is square-integrable as well as continuous. Therefore, HSO is a compact operator. Further, we know that any compact operator possess a singular value decomposition (see \cite{Engl}). Consequently, we may write
 $$
\mathcal{S}(\Psi) = \sum_{k=1}^\infty s_k \langle \Psi, \beta_k \rangle \alpha_k,
$$
where $\{s_k\}$ are singular values of $\mathcal{S}$ satisfying $s_k\to 0$ as $k\to \infty$ and $\{\alpha_k\}, \{\beta_k\}$ are orthonormal systems.
For the operator $\mathcal{S}$, the inverse operator can be written as
 \begin{equation}
 \label{1.3} \Psi = \mathcal{S}^{-1}(\Phi) = \sum_{k=1}^\infty \frac{1}{s_k} \langle \Phi, \alpha_k \rangle \beta_k. 
 \end{equation}
We assume the availability of perturbed data since exact data is not available in practice. So, we may write
$$\Phi = \mathcal{S}(\Psi) + \mathcal{E}.$$ Here $\mathcal{E}$ represents the small error.
Applying $\mathcal{S}^{-1}$ on this $\Phi$ using (\ref{1.3}) to derive that 
\begin{equation}\label{1.4}
\mathcal{S}^{-1}(\Phi) = \Psi+\mathcal{S}^{-1}(\mathcal{E})=\Psi + \sum_{k=1}^\infty \frac{1}{s_k} \langle \mathcal{E}, \alpha_k \rangle \beta_k.
\end{equation}
Since the operator $\mathcal{S}$ is compact, we note that $$s_k\to 0\ \ \ \text{as}\ \ \ \ k\to \infty\implies  \frac{1}{s_k}\to \infty\ \ \ \text{as}\ \ \ \ k\to \infty.$$
This and (\ref{1.4}) derive that $$\|\mathcal{S}^{-1}(\Phi)-\Psi\|_{L^2[0,1]}>C$$ for any positive constant $C>0$. Therefore, there is no continuous dependence between the data and the solution. Hence, the operator $\mathcal{S}$ is ill-posed in the sense of Hadamard due to violation of condition (iii). 
 \subsection{Degree of Ill-posedness}\ 
 In this subsection, we discuss about the operator equations of the form (\ref{1}) with the additional constraint that $T$ is a compact operator between Hilbert spaces. The singular value decomposition (SVD) theorem implies that $T$ has singular values
\[s_1\geq s_2\geq s_3\geq \cdots\geq s_k>0.
\] We note that (\ref{1}) is ill-posed if
$$s_k\to\ 0 \ \ \text{as}\ \ k\to \infty.$$This is also evident from (\ref{1.4}) as decaying values of $s_k$ leads to instability of $A^{-1}$. Further, the degree of ill-posedness can be described by looking at rate at which singular value decays. Accordingly, we have the following two categories of ill-posed problems.
\begin{enumerate}[(i)]
\item \textbf{Mildly ill-posed:} If the decay is of the form
\[
s_k \sim k^{-t}, \quad t > 0,
\]
i.e., $s_k$ decays polynomially (or algebraic decay), then (\ref{1}) is mildly ill-posed.   
\item \textbf{Severely ill-posed:} If the decay is of the form
\[
s_k \sim e^{-r k}, \quad r > 0,
\]i.e., $s_k$ decays exponentially, then (\ref{1}) is severely ill-posed.\end{enumerate}

We note that faster the rate at which $s_k'$s decay, more the degree of ill-posedness of (\ref{1}).
\subsection{Structural Analogy}\ \\
In this subsection, we discuss the analogy between LWE problem and ill-posed inverse problems.
 This is discussed in the following Table \ref{t2}.
 \begin{figure}[!htb]
 \begin{tabular}{c|c|c} 
\textbf{Parameters}& \textbf{LWE Problem} & \textbf{Inverse Problem for HSO}\\ \hline
\textbf{Dimension} & Finite: modular matrix multiplication & Infinite: Operator $\mathcal{S}: L^2[0,1]\to L^2[0,1]$ \\ \hline \hline
 \textbf{Data} &$ \mathbf{b} = \mathbf{A}\mathbf{s} + \mathbf{e}$ with noise $\mathbf{e}$ in finite & $\Phi = \mathcal{S}(\Psi) + \mathcal{E}$ with noise $\mathcal{E}$ in infinite\\  & dimensional space & dimensional space\\ \hline \hline
\textbf{Solution} & Deduce $\mathbf{s}$ from noisy samples $\mathbf{b}$ but& Deduce $\Psi$ from noisy data $\Phi$.  Noise \ \  \ \ \\ &  it is computationally hard due to     & presence makes it highly unstable  and \\  & presence of noise $\mathbf{e}$ & therefore regularization techniques are\\ & &  needed to  reconstruct the exact solution\\ \hline\hline
\textbf{Noise} & Discrete in nature and makes deduction & Analytic and continuous in nature and, \\  & of secret $\mathbf{s}$ computationally infeasible & in practice, noisy data is available. \\ & Here, noise is added deliberately & in place of exact data  \\ \hline \hline 
 \end{tabular}
\captionof{table}{ Analogy between LWE problem and Inverse problems}\label{t2}
\end{figure} 

We note the following:
\begin{itemize}
\item Both LWE problem and Inverse problem for HSO demand to invert a linear operator, which is contaminated by noise. As a result, inversion is either unstable or computationally hard.
\item In inverse problem for HSO, a small noise $\mathcal{E}$ makes the problem highly ill-posed. In LWE problem, noise plays the role of mask, which makes it computationally infeasible to deduce the secret.
\end{itemize}
In this manner, we can see that LWE problem is a structured inverse problem. Precisely, LWE problem can be seen as a special case of inverse problems, i.e., $$\text{LWE Problem}\in \{\text{P}:\ \text{P is an ill-posed inverse problem}\}.$$ In inverse problems, we aim for approximate solution (which is close to exact solution in terms of norm) using regularization techniques whereas for LWE problem, we look for the exact solution. If noise samples in inverse problems are discrete in nature, then LWE problem and inverse problem with matrix operator are same.
\section{Symmetric encryption scheme based on Ill-posed problems}\noindent 
In this section, we propose a symmetric encryption scheme based on ill-posed inverse problems. This is defined through the following subsections.   \\
\textbf{2.1 Parameters and Key Generation}: \begin{itemize}
\item Let $\mathcal{S}: L^2[0,1]\to L^2[0,1]$ be a compact operator. Let $\mathcal{S}^{-1}$ be the inverse of $\mathcal{S}$. We remark that both $\mathcal{S}$ and $\mathcal{S}^{-1}$ are known publically. 
\item Let $\mathcal{E}$ be the error and let it follow certain distribution, e.g., discrete Gaussian distribution or centered binomial distribution. This error is chosen secretly and acts as a secret key. We note that the set containing errors has cardinality atleast $2^{128}$.
\end{itemize}
  \textbf{2.2 Encoding}: \begin{itemize}
  \item   Let $ \{0,1\}^{t}$ be the message space containing all bit-strings of  length $t$. 
  \item Let $\mu\in \{0,1\}^t$ be an arbitrary message string. Then one may see the space $\{0,1\}^t$ as a subset of $L^2[0,1]$ using one of the following two mappings.
  \item \textbf{Map-1}: Using Fourier or Haar basis, generate an orthonormal sequence $\{e_k\}\subset L^2[0,1]$. We enumerate elements in $\{0,1\}^t$ as $\{\mu_1, \mu_2, \cdots, \mu_{2^t}\}$. Then consider the map
  \begin{equation}\label{w1}
    \wp_1: \{0,1\}^t\to L^2[0,1]\ \ \text{defined as}\ \  \wp_1(\mu_k)=e_k. 
\end{equation}
  This map is injective.
  \item \textbf{Map-2}: This map is defined using piecewise constant functions. Consider a string $\mu=s_1s_2\cdots s_t\in \{0,1\}^t$. We define $\Psi_{\mu}\in L^2[0,1]$ as
  $$\Psi_{\mu}(y)=   s_i, \ \qquad \text{if}\ y\in \bigg[ \frac{j-1}{t}, \frac{j}{t}\bigg), \ 1\leq j\leq t. 
  $$
  It can be verified that $\Psi_{\mu}\in L^2[0,1]$.  Accordingly, we define the map
   \begin{equation}\label{w2}
       \wp_2: \{0,1\}^t\to L^2[0,1]\ \ \text{defined as}\ \  \wp_2(\mu)=\Psi_{\mu}.
\end{equation}
  This map is also injective.
  \item Clearly, $\wp_1$,  $\wp_2$ represent the desired encoding of binary messages as elements of $L^2[0,1]$.
\end{itemize}
 \textbf{2.3 Encryption}: 
 \begin{itemize}
 \item
 Let $\mu=s_1s_2\cdots s_t\in \{0,1\}^t$ be the message to be encrypted.  
\item  Apply $\wp_1$ given by (\ref{w1}) (or $\wp_2$ given by (\ref{w2})) on   $\mu$ to obtain $\wp_1(\mu)\in L^2[0,1]$.
\item  The ciphertext is \begin{equation}\label{c}
 \mathcal{C}=\mathcal{S}(\wp_1(\mu))+\mathcal{E}, 
\end{equation} where $\mathcal{E}$ is the secret key. This error introduces noise in the data. 
 \end{itemize}
  \textbf{2.4 Decryption and Decoding}: \begin{itemize}
  \item The decryptor after receiving  $\mathcal{C}$, uses secret key $\mathcal{E}$ to obtain
  $$\mathcal{C}-\mathcal{E}=\mathcal{S}(\wp_1(\mu))$$
  \item Apply $\mathcal{S}^{-1}$ on the exact data to obtain 
  $$\mathcal{S}^{-1}(\mathcal{C}-\mathcal{E})=\wp_1(\mu).$$
  \item Represent $\wp_1(\mu)$ as an element of $\{0,1\}^t$ by using $\wp_1^{-1}$ (which exists on its range). This yields the message.
  \end{itemize}
  Next, we discuss the security of our encryption scheme along with   its certain characteristics.\\
   \textbf{2.5 Security Analysis}:\\
  \textbf{(a) Brute force attack}:\\ The adversary needs to try all the possible errors $\mathcal{E}$ from the set of errors to get the message from the knowledge of $\mathcal{C}$. This is computationally infeasible to perform in polynomial time. \\
  \textbf{(b) Error should be ephemeral}:\\ Suppose two messages $\mu_1$ and $\mu_2$ are such that same error $\mathcal{E}$ is taken to encrypt both. Then, (\ref{c}) implies that
   $$\mathcal{C}_1=\mathcal{S}(\wp_1(\mu_1))+\mathcal{E}$$
  $$ \mathcal{C}_2=\mathcal{S}(\wp_1(\mu_2))+\mathcal{E}.$$
  These two imply that
  $$ \mathcal{C}_1-\mathcal{C}_2=\mathcal{S}(\wp_1(\mu_1))-\mathcal{S}(\wp_1(\mu_2)).$$
  If the operator $\mathcal{S}$ is linear, then this means 
  \begin{equation}\label{e1}
     \mathcal{C}_1-\mathcal{C}_2=\mathcal{S}(\wp_1(\mu_1)-\wp_1(\mu_2)). 
\end{equation}
  Further, if $\wp_1$ is also linear, then (\ref{e1}) gives the encryption of $\mu_1-\mu_2$. Therefore, for every message, error $\mathcal{E}$ should be chosen uniformly at random.\\
   \textbf{(c) Probabilistic encryption}:\\ Our scheme comes under the category of probabilistic encryption. Given a fixed message $\mu$, there can be many ciphertexts corresponding to this message. This is because by changing the error, ciphertext gets changed corresponding to a fixed message.\\
 \textbf{(d) Security against ciphertext only attack}:\\
 The adversary $\mathfrak{A}$ knows $\mathcal{C}$ and $\mathcal{S}^{-1}$. Using these in (\ref{c}), $\mathfrak{A}$ computes
 $$ \mathcal{S}^{-1}(\mathcal{C})=\wp_1(\mu)+\mathcal{S}^{-1}(\mathcal{E}).$$
 But due to ill-posedness of the compact operator equation $\mathcal{S}(\Psi)=\Phi$ (or inverse problem), the small error $\mathcal{E}$ leads to large error in $ \mathcal{S}^{-1}(\mathcal{C})$ as shown in subsection 1.2. Therefore, it becomes computationally infeasible for $\mathfrak{A}$ to deduce $\mu$ without applying regularization techniques. We remark that the application of regularization techniques further depends on the degree of ill-posedness of the inverse problem as discussed in subsection 1.3.\\
 \textbf{(e) Security against known/chosen plaintext attacks}:\\
  The adversary knows (or demands) polynomially many plaintext-ciphertext pairs $(\mu_k, \mathcal{C}_k)$ (these may be chosen adaptively). The main purpose of $\mathfrak{A}$ is to derive $\mathcal{E}$. It follows from (\ref{c}) that
  $$ \mathcal{C}_k=\mathcal{S}(\wp_1(\mu_k))+\mathcal{E}_k, \ \ k\geq 1.$$
  Since $\mathcal{E}_k$ is different and randomly chosen for every message, it becomes computationally infeasible to deduce $\mathcal{E}$ for a new message.\\
  \textbf{(f) Security against Differential cryptanalysis}:\\
  Since we have shown that our scheme is safe against adaptive chosen plaintext attack, therefore, our scheme is also resistant to Differential cryptanalysis.\\ \textbf{(g) CCA2 security (security against adaptive chosen ciphertext attack)}:\\
  Due to inclusion of unique/random errors at each instance, it would not be possible for  $\mathfrak{A}$ to deduce the key/encryption of new message from previous adaptive plaintext-ciphertext queries in polynomial time. Thus, the scheme is CCA2 secure.
 \\
 \textbf{Remark 1} (Requirement of a synchronized error generator): For a symmetric encryption scheme, both sender and receiver should have the same secret key. In our case, the secret key is error. Consequently, both the parties should have  a synchronized  error generator that generates error uniformly at random (or following certain distribution).  
  \section{Public key encryption scheme based on Ill-posed problems}\noindent In this section, we build on the work of previous section and propose a public key encryption (PKE) scheme based on ill-posed problems. Specifically, we integrate the symmetric encryption scheme proposed in the previous section with CRYSTALS-Kyber key encapsulation mechanism (KEM) \cite{Kyber}.  The PKE is defined as follows (see  also Table \ref{t1}). \\
  \textbf{3.1 Key Generation}:\\
  We run CRYSTALS-Kyber KEM
   \emph{key generation algorithm} to generate a public/private key pair. We denote it by pk and sk, respectively.\\
    \textbf{3.2 Encryption}:\\
    Let $\mu\in \{0,1\}^t$ be the message to be encrypted. We run CRYSTALS-Kyber KEM \emph{encaps algorithm} to generate a key $\mathcal{K}'$  (of length $256$ bits) and  the corresponding ciphertext $\mathcal{C}_1$. Further, we instantiate XOF function \cite{fips202} on $\mathcal{K}'$ to get a key $\mathcal{K}$ of desired length to be used for symmetric encryption scheme $\mathsf{S}$ described in section 2.   Finally, we encrypt $\mu$ using the encryption function of the scheme $\mathsf{S}$ to get ciphertext $\mathcal{C}_2$. The final ciphertext is $(\mathcal{C}_1, \mathcal{C}_2)$.\\
    \textbf{3.3 Decryption}:\\
    After receiving $(\mathcal{C}_1, \mathcal{C}_2)$, the decryptor first run CRYSTALS-Kyber KEM \emph{decaps algorithm} with input $\mathcal{C}_1$ to get the key $\mathcal{K}'$. Then, the decryptor runs  XOF function   on $\mathcal{K}'$ to get the key $\mathcal{K}$.  Finally,   using  the decryption function of the scheme $\mathsf{S}$ on ciphertext $\mathcal{C}_2$ with key $\mathcal{K}$, we  get the message $\mu$.
   
\begin{figure}[!htb]
	\scalebox{0.9}{	\centering
		\begin{tabular}{|c|c|c|}
			\hline\hline			
		\textbf{Key Generation}&	\textbf{Encryption}  & \textbf{Decryption}    \\ \hline \hline 
	 		$\bullet$ Run CRYSTALS-Kyber key & &\\  generation algorithm & &\\	
	 		$\bullet$ Output: (pk, sk)\ \ \qquad\qquad\qquad & $\bullet$ pk is known to encryptor  \qquad\quad & $\bullet$ sk is known to decryptor\ \qquad\ \ \  \\ \hline
			& $\bullet$ Message $\mu\in \{0, 1\}^t$\ \ \qquad \qquad \qquad &\\
			& $\bullet$ Run CRYSTALS-Kyber enca- &\\ & ps algo with pk to get $(\mathcal{K}', \mathcal{C}_1)$ &\\
 & $\bullet$ Apply XOF on $\mathcal{K}'$ to get $\mathcal{K}$\ \ \  \ \  &\\
 & $\bullet$ Encrypt $\mu$  through     $\mathsf{S}$ via key \ \    &\\ &  $\mathcal{K}	$ to get ciphertext $\mathcal{C}_2$ 	&\\	& $\bullet$ The ciphertext is $(\mathcal{C}_1, \mathcal{C}_2)$\ \qquad\qquad &\\ \hline
 & & $\bullet$ Received ciphertext $(\mathcal{C}_1, \mathcal{C}_2)$ \  \ \  \  \     \ \\
  & &   $\bullet$ Run CRYSTALS-Kyber deca-  \ \ \  \\ & & ps algo on $\mathcal{C}_1$ with sk  to  get  $\mathcal{K}'$   \\ 
  & & $\bullet$ Apply XOF on $\mathcal{K}'$ to get $\mathcal{K}$\ \  \  \  \ \  \ \\ & &    $\bullet$ Deccrypt $\mathcal{C}_2$  through    $\mathsf{S}$ via       $\mathcal{K}	$   \  \ \  \  \\
			\hline
		\end{tabular}	}
	\captionof{table}{ Complete description of our PKE scheme}\label{t1}
	\end{figure}\noindent    \textbf{3.4 Security Analysis}:\\
    The security of our PKE clearly depends on the security of Kyber KEM and symmetric encryption scheme $\mathsf{S}$. The security of Kyber KEM depends on  lattice problem, i.e., module-LWE (module learning with errors), whose hardness is very well understood (see \cite{Kyber}). The security of $\mathsf{S}$ is already discussed, i.e., it is CCA2 secure. Consequently, it follows from \cite{Kyber} that our PKE is IND-CCA2 secure (indistinguishability under adaptive chosen ciphertext attack), which is the golden security standard.    
 \section{Discussion}\noindent
We have unveiled a  canonical analogy between post-quantum lattice based learning with error problem and ill-posed inverse problems. Precisely, we have shown that LWE problem is a special case of solving inverse problems. Motivated from this fact, we have proposed two encryption schemes. The first one is symmetric and other one is aymmetric (or PKE). We   also thoroughly discussed the security of these schemes. 
In future, this work can be extended in a number of ways. The first one is to look at the impact of regularization techniques on the security of these schemes. The second one is to look at the efficient ways of sampling errors in (\ref{c}).

\bibliographystyle{plain}

\end{document}